\documentclass[preprint,10pt,times,natbib,fleqn,runningheads,authoryear]{article}
\usepackage{titlesec}
\titleformat{\section}
        {\normalfont\filright\normalsize} 
        {\thesection}
        {0.8em}
        {\bf}[]
        
    \titleformat{\subsection}
        {\normalfont\filright\normalsize} 
        {\thesubsection}
        {0.8em}
        {\bf}[]

\usepackage{bm}

\usepackage{amsmath}
\usepackage{amssymb}
\usepackage{theorem}
\usepackage{xspace}
\usepackage{geometry}
\usepackage{array}
\usepackage{subfigure}
\usepackage{epsfig}
\usepackage{hhline}
\usepackage{bbm}
\usepackage{natbib}   
\usepackage{setspace}
\usepackage[utf8]{inputenc}    
\DeclareMathAccent{\widehat}{\mathord}{largesymbols}{"62}
\DeclareMathAccent{\widetilde}{\mathord}{largesymbols}{"65}
\def\pth#1{\left(#1\right)}

\def\cro#1{\left[#1\right]}

\def\eeX{\mathbb{X}}

\numberwithin{equation}{section}   

\def\ebo{\textrm{\mathversion{bold}$\mathbf{\beta^0}$\mathversion{normal}}}
\def\uu{\textrm{\mathversion{bold}$\mathbf{1}$\mathversion{normal}}}

\def\eb{\textrm{\mathversion{bold}$\mathbf{\beta}$\mathversion{normal}}}  
\def\ed{\textrm{\mathversion{bold}$\mathbf{\delta}$\mathversion{normal}}}

\def\eU{\textrm{\mathversion{bold}$\mathbf{\Upsilon}$\mathversion{normal}}} 

\def\eE{I\!\!E}

\def\e1{1\!\!1}

\def\XX{\textrm{\mathversion{bold}$\mathbf{X}$\mathversion{normal}}}

\def\xx{\textrm{\mathversion{bold}$\mathbf{x}$\mathversion{normal}}}

\theoremstyle{plain}

\newtheorem{theorem}{Theorem}[section]

\newtheorem{lemma}{Lemma}[section]

\newtheorem{remark}{Remark}[section]

\newcommand{\beqn}{\begin{eqnarray*}}
\newcommand{\eeqn}{\end{eqnarray*}}
  
\def\ee1{\textrm{\mathversion{bold}$\mathbf{\varepsilon}$\mathversion{normal}}}

\newcommand{\bfz}{{\bf z}}

\def\eu{\mathbf{{u}}}

\newcommand{\N}{\mathbb{N}}
\newcommand{\R}{\mathbb{R}}
\newcommand{\PP}{\mathbb{P}}

\def\argmin{\mathop{\mathrm{arg\,min}}} 
\def\hh{ \hspace*{0.5cm}}
\usepackage{fancyhdr}   
 
\begin{document}
{\large \title {{\bf Adaptive Fused LASSO in  Grouped Quantile Regression}}

\author{GABRIELA CIUPERCA  \footnote{Universit\'e Claude Bernanrd Lyon 1,  Institut Camille Jordan, 
Bat.  Braconnier, 43, blvd du 11 novembre 1918, 
F - 69622 Villeurbanne Cedex, France,  E-mail: {\it Gabriela.Ciuperca@univ-lyon1.fr},}\\
Universit\'e Claude Bernanrd Lyon 1,  UMR 5208, Institut Camille Jordan,  France
}}
\maketitle
\noindent {\bf \hrule }
 \begin{abstract}
{\large This paper considers quantile model with grouped explanatory variables. In order to have the sparsity of the parameter groups but also the sparsity between two successive groups of variables, we propose and study an adaptive fused group LASSO quantile estimator. The number of variable groups can be fixed or divergent. We find the convergence rate   under classical assumptions and we show that the proposed estimator satisfies the oracle properties. }
\end{abstract}
\noindent {\bf Keywords:} {\large
 group selection; quantile regression; adaptive fused LASSO; selection consistency; oracle properties. 
}
\noindent {\bf \hrule }

\maketitle
\newpage
\section{Introduction}
\hh The idea of this paper comes from the ascertainment that in many practical applications, for studying a process or a random variable in function of grouped explanatory variables, we want to identify significant groups of variables but also to make a hierarchy between these groups. The explanatory variables can be continuous or discrete.  The most common example of linear model with grouped variables is the multivariate variance analysis. But in many situations for theoretical study of the linear models, classical assumptions are imposed on errors: zero mean and bounded variance, which is not often the case in applications. Then, if these classical assumptions   are not satisfied or if the model has heavy-tailed errors, a very interesting approach is the quantile method. Moreover, compared to classical estimation methods (least squares, least absolute deviations) which give the model behaviour around of the mean or of the median, the quantile method offers a very complex and global   insight. This method allows to study how the explanatory variables influence the response variable distribution.   For a complete overview on quantile method, we refer the reader to book of \cite{Koenker-05}. 
 Then, in order to cover more possible cases of models, in this paper, we will consider grouped quantile regression, which allows the relaxation of the classical conditions on the two first moments of the model error.  We also want to identify the relevant variable groups, an automatic selection by a LASSO type method of the significant variable groups being more meaningful that an automatic selection of individual variables.  The LASSO estimator, introduced by \cite{Tibshirani-96} for the least squares framework, doesn't always satisfy the automatic selection, and then a solution is the  adaptive LASSO estimator, proposed initially by \cite{Zou-06}. \\
The LASSO methods have been the subject of active research in the last decade. We give here only the references concerning the LASSO methods  for models with grouped variables. 
\cite{Zhang-Xiang-15}, \cite{Ciuperca-16} have considered the adaptive group LASSO in high-dimensional linear model by penalizing the sum of squares, respectively quantile process. Earlier, \cite{Wei-Huang-10} had considered the  adaptive group LASSO estimator but for gaussian model error.  \cite{Wang-You-Lian-15} study also the convergence and the sparsity of the (non adaptive) group LASSO estimator in a high-dimensional generalized linear model. When the number of groups is fixed, \cite{Wang-Leng-08} show the model consistency obtained by an adaptive group LASSO method. For a review, but until 2012, of group selection methods and several applications of these methods the reader can see \cite{Huang-Breheny-Ma-12}.\\
\hh For a linear model, without grouped variables, in order to encourage sparsity of the parameters but also the sparsity of their differences, to identify predictive variable clusters,   \cite{Tibshirani-Saunders-05} introduced an additional penalty to the LASSO penalty, taking the $L_1$ norm of the differences between two successive  parameters. They called the obtained estimator, fused LASSO estimator. Applications of the  fused LASSO method  can be found in  \cite{Jang-Lim-Lazar-15} or in \cite{Li-Zhu-07}. The fused LASSO idea was adopted by  \cite{Jiang-Wang-Bondell-13}, \cite{Jiang-Bondell-Wang-14}, \cite{Zhao-Zhang-Liu-14} to quantile model.   A very recent paper of \cite{Viallon-Lambert-Picard-16}, considers a generalized linear models, estimated by adaptive fused LASSO method. Adaptive fused LASSO penalty is also used by \cite{Sun-Wang-Fuentes-16} for estimation of the spatial and temporal quantile functions. These papers have in common that the models have fixed  number of explanatory variables.\\
\hh The fused LASSO penalty for grouped variables, proposed and studied in this paper,  will be helpful to strengthen the sparsity between two successive groups of variables.  To the knowledge of the author, the adaptive fused LASSO method wasn't considered for a linear model and further, with the possibility that the number of groups converges to infinity when number of observations diverges. Even for quantile linear models without grouped variables, there is no work in literature on the adaptive fused LASSO method.  This is the originality of the present paper. Emphasize that the proposed estimator and obtained results are valid for a wide spectrum of error distributions.\\

The paper is organized as follows. In Section \ref{section 2}, we present the quantile model with grouped  variables,  we also introduce the adaptive fused estimator and we give general notations and assumptions. In Section \ref{section 3}, we study the convergence rate, oracle properties of the estimator when the group number is fixed. A general convergence rate when the group number diverges is found in Section \ref{section 4}. In the same section, we state that the oracle properties remain true. All proofs will be postponed in Section \ref{preuves}.
\section{Model, notations and assumptions}
\label{section 2}
\hh In this section, we first introduce the quantile model with grouped variables. Afterwards, some notations used throughout in the paper are given, followed by the introduction of the proposed fused  estimator. Finally, general assumptions on the model errors, design and on the group number are given. \\ 
\hh Let us consider the following linear model with $p$ groups of variables
\begin{equation}
\label{eq1}
Y_i=\sum^p_{j=1}\XX_{ij}^t \eb_j +\varepsilon_i=\eeX_i^t \eb +\varepsilon_i, \qquad  i =1, \cdots, n. 
\end{equation}
The random variables in model (\ref{eq1}) are: $Y_i$ the response variable and  $\varepsilon_i$ the model error. The column vector $\eeX_i$ is the $i$th observation of the explanatory variables and its contains $p$ groups of variables.  For each group  $j$,  with  $j =1, \cdots, p$, the vector of the parameters is  $\eb_j \equiv (\beta_{j1}, \cdots , \beta_{j d_j}) \in \R^{d_j}$ and the design $\XX_{ij}$ for observation $i$,   is a vector of size $d_j \times 1$. The vector with the all  coefficients is  $\eb\equiv (\eb_1, \cdots, \eb_p) $  and  $ \eb^0_j=(\beta^0_{j1}, \cdots ,\beta ^0_{j d_j})$ the true (unknown) value of the parameter $\eb_j$, for $j =1, \cdots, p$.  For observation $i$, we denote by $X_{ij,k}$ the  $k$th variable of the $j$th group. We will assume that $d_j=d$ for any $j=1, \cdots, p$ by taking $d=\max_{j=1, \cdots , p} d_j$, filling the  components of $\eb_j$ between $d_j$ and $d$ with  0 and the values of $\eeX_i$ also with 0.  \\
For model (\ref{eq1})  it is possible that there are insignificant variables groups. For this, we will consider the index set of the significant groups: 
  $${\cal A} \equiv \{ j \in \{ 1, \cdots p\} ; \| \eb^0_j \| \neq 0\}$$
   and obviously the index set of the insignificant groups ${\cal A}^c \equiv \{ j ; \| \eb^0_j \| = 0\}$. We denoted by $|{\cal A}|$ the cardinal of the index set ${\cal A}$.  Obviously, in practical applications, the two sets  ${\cal A}$ and ${\cal A}^c$ are unknown.   \\
On the other hand, we denote by $|{\cal A}|=p^0$, $r^0=d p^0$ and $r=dp$. The numbers $p$ and $d$ are known,  $p^0$ is contrariwise unknown and then  $r^0$ also.   \\
For a $r$-vector of parameters $\eb$, we denote $\eb_{\cal A}$ the subvector of $\eb$, of dimension $r^0 \times 1$, which contains $\eb_j$, for $j\in {\cal A}$. The $(r-r^0)$-vector $\eb_{{\cal A}^c}$   contains  $\eb_j$ for $j\in {\cal A}^c$.\\

We introduce now the quantile model. This method allows the non necessity of  the classical assumptions on errors: $\eE[\varepsilon_i]=0$ and $Var(\varepsilon_i) < \infty$. Since these assumptions are not often satisfied in practical applications,  the quantile method can be used  extensively   in many different areas.\\
So, for a  quantile index $\tau \in (0,1)$, the check function $\rho_\tau (.): \R \rightarrow \R_+$ is defined by $\rho_\tau(u)=u(\tau - \e1_{u <0})$. In this paper, the index $\tau$ is considered fixed.\\
 
 Before defining the adapted fused LASSO estimator for the parameter $\eb$ of (\ref{eq1}), we give some general notations. All throughout the paper,  $C$ denotes a positive generic constant not dependent on $n$, which may take different values in different formula or even in different parts of the same formula. The value of $C$ is not of interest. All vectors and matrices are denoted by bold symbols and all vectors are written as column vectors. For a   vector $\textbf{v}$, we denote by $\textbf{v}^t$ its transposed and by $\|\textbf{v} \|$ its Euclidean norm. Notations $ \overset{\cal L} {\underset{n \rightarrow \infty}{\longrightarrow}}$, $ \overset{\PP} {\underset{n \rightarrow \infty}{\longrightarrow}}$ represent the convergence in distribution and in probability, respectively, as $n \rightarrow \infty$.  For a positive definite matrix $\textbf{M}$, we denote by $\lambda_{\min}(\textbf{M)}$  and $\lambda_{\max}(\textbf{M)}$ its the smallest and largest eigenvalues, respectively. When it is not specified, the convergence is for $n \rightarrow \infty$.\\
 
For model (\ref{eq1}), $n$ observations of $(Y_i, \eeX_i)_{1 \leqslant i \leqslant n}$ are available.  In order to define the estimator that will allow automatic selection of significant groups of variables,  we must first consider the quantile process:
\[
G_n(\eb) \equiv \sum^n_{i=1} \rho_\tau(Y_i-\eeX^t_i \eb).
\]
The quantile estimator for  $\eb$ is the minimizer of the  quantile process:
  \begin{equation}
  \label{eq2}
  \widetilde{\eb}_n \equiv \argmin_{\eb \in \R^r} G_n(\eb).
  \end{equation}
 This estimator can be written taking into account each group $\widetilde{\eb}_n=(\widetilde{\eb}_{n;1}, \widetilde{\eb}_{n;2}, \cdots , \widetilde{\eb}_{n;p})$, with $\widetilde{\eb}_{n;j}$ a  vector of size $d$, for $j =1, \cdots, p$. We will use $\widetilde{\eb}_n$ for constructing the two adaptive LASSO penalties. Note that, the number $r$ of the total variables needs to be smaller than the sample size $n$.  For model (\ref{eq1}), we define the adaptive fused group  LASSO quantile (\textit{afg\_LASSO\_Q}) estimator, denoted by  $\widehat{\eb}_n$,  as the minimizer of the following process:
  \begin{equation}
  \label{eq3}
 Q_n(\eb) \equiv G_n(\eb) +\mu_n^{(1)}\sum^p_{j=1}  \widehat{\omega}_{n;j}^{(1)} \| \eb_j\|+\mu_n^{(2)}\sum^p_{j=2}  \widehat{\omega}_{n;j}^{(2)} \| \eb_j-\eb_{j-1}\|,
  \end{equation}
 with the weights  $\widehat{\omega}_{n;j}^{(1)} \equiv \| \widetilde{\eb}_{n;j} \| ^{- \gamma}$, $\widehat{\omega}_{n;j}^{(2)} \equiv \| \widetilde{\eb}_{n;j} -\widetilde{\eb}_{n;j-1} \| ^{- \gamma}$ and  $\gamma >0$ a fixed known parameter.  The estimator  $\widehat{\eb}_n$ is written $\widehat{\eb}_n =(\widehat{\eb}_{n;1}, \cdots, \widehat{\eb}_{n;p})$ and $\widehat{\eb}_{n;j}$ is a vector of size $d$, for $j=1, \cdots , p$. The tuning parameters $\mu_n^{(1)}$, $\mu_n^{(2)}$ are assumed to converge to infinity as $n \rightarrow \infty$. Additional conditions on $\mu_n^{(1)}$, $\mu_n^{(2)}$, taking into account $\gamma$ and the group number  $p$, will be given later.\\
 
  The purpose of this paper is to study the properties of the estimator $\widehat{\eb}_n$, mainly the oracle properties: the significant groups of  variables are estimated, with an optimal estimation rate, by asymptotically gaussian estimators  and  the corresponding parameters to nonsignificant groups are shrunk directly to 0 with a probability converging to one.  
In order to study the asymptotic properties of the   \textit{afg\_LASSO\_Q} estimator $\widehat{\eb}_n$, for some  $r$-vector  $\eb \in \R^r$, we also consider the process:
  \begin{equation}
  \label{Lnb}
  L_n(\eb) \equiv Q_n(\eb)- Q_n(\ebo).
  \end{equation}
 Let us note that, for $\tau=1/2$,   model  (\ref{eq1}) becomes median regression with grouped variables. The estimator  $\widehat{\eb}_n$ becomes in this case, adapted fused grouped LASSO median estimator. \\ 
  
The asymptotic properties for $\widehat{\eb}_n$ will be studied under the following assumptions for errors,  design and group number $p$:\\
\textbf{(A1)}  $(\varepsilon_i)_{1 \leq i \leq n}$ are i.i.d., with  the distribution function  $F$ and density function $f$. The density function $f$ is continuously,  strictly positive in a neighborhood of zero and has a bounded first derivative in the neighborhood of 0. The $\tau$th quantile of $\varepsilon_i$ is zero: $\tau= F(0)$. \\
\textbf{(A2)} There exist constants $0 < m_0 \leq M_0 < \infty$ such that 
 \[ m_0 \leq \lambda_{\min} ( {n}^{-1} \sum^n_{i=1} \eeX_i \eeX_i^t) \leq \lambda_{\max} ( {n}^{-1} \sum^n_{i=1} \eeX_i \eeX_i^t) \leq M_0 . \]
\textbf{(A3)}  $\pth{ {p}/{n}}^{1/2} \max_{1 \leqslant i \leqslant n} \| \eeX_i\| \rightarrow 0$.\\
\textbf{(A4)} $p$ is such that $p=O(n^c)$, with $0 \leq c < 1$.\\
For the smallest nonzero vector norm and on constant $c$ of assumption (A4) we assume:\\
\textbf{(A5)} Let us denote $h_0 \equiv \min_{1 \leqslant j \leqslant p_0} \| \eb^0_j\|$. There  exists a  constant $M>0$ such that $M \leq n^{- \alpha} h_0$ and $\alpha > (c-1)/{2}$. \\
Concerning the size of the nonzero  parameter vectors, we take the following assumption: \\
\textbf{(A6)} $r^0=O(p_0)$.\\

Assumptions (A2), (A3) are standard for LASSO methods and (A1) is classic for quantile regression (see \cite{Ciuperca-15b},  \cite{Koenker-05}, \cite{Zou-Yuan-08}, \cite{Wu-Liu-09}).  
Assumptions (A3), (A4) are also considered in \cite{Ciuperca-16}, \cite{Zou-Zhang-09} for high-dimensional linear model, while (A5) and (A6) are required for adaptive group LASSO least square estimator in \cite{Zhang-Xiang-15} and in \cite{Ciuperca-16} for adaptive group LASSO quantile estimator. \\
For the  case $p$ fixed, then $c=0$, only assumptions (A1) - (A3) will be needed. For the case $p=p_n \rightarrow \infty$ as $n \rightarrow \infty$,  assumptions (A4), (A5) and (A6) are also considered, with $c \in (0,1)$ in assumption (A3). \\

 In Sections \ref{section 3} and \ref{section 4}, we will study the adaptive fused group LASSO quantile (\textit{afg\_LASSO\_Q}) estimator $\widehat{\eb}_n$ for two cases of the group number: $p$ fixed and  $p \rightarrow \infty$ as $n \rightarrow \infty$, respectively. 
  \section{Case $c=0$}
  \label{section 3}
 \hh In this section we will propose and study the asymptotic properties of the   \textit{afg\_LASSO\_Q} estimator of the parameter $\eb$ for model (\ref{eq1}), when the number of groups $p$ is fixed. \\

 Regarding   assumptions, as specified above, in order to prove the oracle properties for $\widehat{\eb}_n$, only (A1), (A2), (A3) will be needed, with a weaker condition in (A1) on error density  $f$. So, the condition that $f$ has a bounded derivative in the neighbourhood of 0 with an weaker condition can be replaced in assumption (A1) by: for every $e \in int({\cal B})$,  $\uu_r \in \R^r$, we have 
 \begin{equation}
 \label{rA2}
 \lim_{n \rightarrow \infty} n^{-1} \sum^n_{i=1} \int^{\xx^t_i \uu_r}_{0} \sqrt{n}[F(e+n^{-1/2}v)-F(e)] dv = \frac{1}{2} f(e) \uu^t_r 
 \eU
 \uu_r .
 \end{equation}
 The $r$-vector  $\uu_r$ contains as elements 1. The matrix $\eU$ is defined by (\ref{MU}). \\
 Note also that assumption (A3) becomes: $n^{-1} \max_{1 \leq i \leq n} \eeX_i^t \eeX_i {\underset{n \rightarrow \infty}{\longrightarrow}}  0$ and assumption (A2) implies that 
  \begin{equation}
  \label{MU}
  n^{-1} \sum^{n}_{i=1} \eeX_i \eeX_i^t {\underset{n \rightarrow \infty}{\longrightarrow}} \eU,\end{equation}
   with $\eU$ a positive definite matrix.\\
   
   The tuning parameters $\mu_n^{(1)}$, $\mu_n^{(2)}$ and the positive  constant $\gamma$ are such that, for $n \rightarrow \infty$, 
\begin{equation}
\label{cond_lambda}
\mu_n^{(m)} \rightarrow\infty,\quad   n^{-1/2}\mu_n^{(m)}  \rightarrow 0, \quad n^{(\gamma-1)/2} \mu_n^{(m)} \rightarrow \infty, \textrm{ for } m=1,2.
\end{equation}
 
 For $m=1$, we get the conditions imposed on the tuning parameter  by \cite{Ciuperca-16} for adaptive group LASSO quantile estimator (non fused). Conditions in  (\ref{cond_lambda}) on $\mu_n^{(1)}$ and $\mu_n^{(2)}$ are also found in \cite{Viallon-Lambert-Picard-16}, where an adaptive fused LASSO for generalized linear models is considered. For the particular case $\gamma=1$, for a  quantile model without grouped variables ($d=1$), we obtain the conditions on the tuning parameters imposed by  \cite{Jiang-Bondell-Wang-14}.\\
   
The proofs of all results are given in Section \ref{preuves},   sub-section \ref{subsection 5.1}.\\ 
 By the following lemma we show that, when the variables are grouped, the adapted fused group LASSO quantile parameter estimator has the same convergence rate   as by classical quantile method, without grouping variables, without adapted fused LASSO penalty.  This convergence rate will serve as an essential tool for studying process $L_n(\eb)$ when  $\eb$ belongs to a neighbourhood of $\ebo$ of order radius $n^{-1/2}$ and for showing the asymptotic normality of the parameter estimators corresponding to the significant groups.
\begin{lemma}
\label{Lemma 2.1}
Under assumptions  (A1), (A2), (A3) and conditions in (\ref{cond_lambda}) for the  tuning parameters, we have, $n^{1/2} \|\widehat{\eb}_n- \eb^0\|=O_{\PP}\pth{1}$.
\end{lemma}

In order to study the oracle properties of the estimator $\widehat{\eb}_n$, let us consider the index set of the groups selected by the following adaptive fused group  LASSO quantile method:
\[
\widehat{\cal A}_n \equiv \{ j \in \{1, \cdots, p\}; \; \| \widehat{\eb}_{n;j}\| \neq 0 \}  
\]
and ${\widehat{\cal A}_n}^c$ its complementary.\\

The following theorem shows a first oracle property, that the \textit{afg\_LASSO\_Q} estimators with indices in the set ${\cal A}$ are asymptotically  Gaussians. 

\begin{theorem}
\label{theorem 1}
Under assumptions (A1), (A2) and (A3) and conditions of (\ref{cond_lambda}),  we have $\sqrt{n}( \widehat{\eb}_n - \eb^0)_{\cal{A}}  \overset{\cal L} {\underset{n \rightarrow \infty}{\longrightarrow}} {\cal N}\big(\textbf{0}_{r^0}, \tau(1-\tau) f^{-2}(0)\eU^{-1}_{\cal A} \big)$, with $\eU_{\cal A}$ the submatrix of $\eU$ with the row and column indices in ${\cal A} $. 
\end{theorem}

Compared to the adaptive LASSO quantile method, if an additional penalty fused is considered, we got the same variance matrix for the asymptotic gaussian law  (see \cite{Ciuperca-16}). \\
In practical applications, the set ${\cal A}$ is unknown. In exchange, it can be estimated by $\widehat{\cal A}_n$. Then, for that the estimator $\widehat{\eb}_n$ to be interesting, it is necessary that these two sets coincide with probability converging to 1, as $n$ converges to infinity. 
By the following theorem, we show that the second oracle property, i.e. the sparsity, is satisfied for the  \textit{afg\_LASSO\_Q} estimator.

\begin{theorem}
\label{th_selection}
Under the same  assumptions as in Theorem \ref{theorem 1}, we have,  $\lim_{n \rightarrow \infty}\PP [\widehat{\cal A}_n={\cal A}] =  1 $.
\end{theorem}

The proof of Theorem \ref{th_selection}, given in sub-section \ref{subsection 5.1},   is in two parts. The result $\lim_{n \rightarrow \infty}\PP [{\cal A} \subseteq \widehat{{\cal A}}_n] =  1$ is an immediate consequence of Theorem \ref{theorem 1}. In order to prove $\lim_{n \rightarrow \infty}\PP [Card(  {\cal A}^c \cap \widehat{{\cal A}}_n) \geq 1]=0 $, the proof is  quite technical, taking into account already proven properties to $\widehat{\eb}_n$ and imposed conditions for the tuning parameters.

\section{Case $c>0$}
\label{section 4}
\hh In this section we consider same model (\ref{eq1}) with grouped  variables, but with the number  $p$ of groups depending on $n$ and divergent: $p=p_n$ and $p_n \rightarrow \infty$ as $n \rightarrow \infty$.  For readability we keep notation $p$ instead of $p_n$. Similarly for $r=p d$. The main purpose is to show that the \textit{afg\_LASSO\_Q} estimator keeps the oracle properties even though the group number diverges.  The proofs of all results are given in Section \ref{preuves},  sub-section \ref{subsection 5.2}. A major difficulty that appears in the proofs is that the size of vectors and of matrices converges to infinity when $n$ tends to infinity.\\

In order to show the main result of this Section, we will first find the  convergence rate of adaptive fused group  LASSO quantile estimator $\widehat{\eb_n} $ of $\eb$. Afterwards, we will show that this estimator satisfies the oracle properties. We recall that the two tuning parameters $\mu^{(1)}_n$ and $\mu^{(2)}_n$  converge to infinity as $n \rightarrow \infty$. 

\begin{lemma}
\label{th_vconv}
Under assumptions (A1)-(A5)  and the two tuning parameters $(\mu_n^{(m)})_{n \in \N}$ satisfying $\mu^{(m)}_n n^{(c-1)/{2} - \alpha \gamma} \rightarrow 0$, as $n \rightarrow \infty$, for $m=1,2$, we have   $\|\widehat{\eb}_n-\eb^0\|=O_{\PP}\pth{ ( {p}n^{-1})^{1/2}}$.
\end{lemma}

We observe that for fixed $p$, we obtain the result of  Lemma \ref{Lemma 2.1}. The convergence rate  as $p \rightarrow \infty$ of the \textit{afg\_LASSO\_Q} estimator is the same as  that of  \cite{Ciuperca-16} for adaptive group LASSO quantile estimator. Then, the fused penalty doesn't affect the estimator rate convergence.  For the particular case $c=\alpha=0$, the condition imposed on $\mu^{(m)}_n$, for $ m=1,2$,   in Lemma \ref{th_vconv}, is the second condition of (\ref{cond_lambda}).\\

In order to prove the sparsity property, the assumptions used in Lemma \ref{th_vconv} are sufficient. 
Since $p \rightarrow \infty$, we need in addition assumption (A6) for showing the asymptotic normality of the    \textit{afg$\_$Q$\_$LASSO} estimators for the significant groups of variables. For the  tuning parameters, we consider a generalization for the third condition of (\ref{cond_lambda}). 

\begin{theorem}
\label{Theorem 2SPL} Suppose that assumptions (A1)-(A5)  are  satisfied  and also that the tuning parameters satisfy   $\mu^{(m)}_n n^{(c-1)/{2} - \alpha \gamma} \rightarrow 0$, $\mu^{(m)}_n n^{\big(-c(1+\gamma)+\gamma-1 \big)/2} \rightarrow \infty$, as $n \rightarrow \infty$, for $m=1,2$.  Then:\\
(i) $\PP \cro{\widehat{\cal A}_n={\cal A}}\rightarrow 1$,  as $n \rightarrow \infty $.\\
(ii) If moreover assumption (A6) holds, for any  vector $\eu$ of size $r^0$ such that $\| \eu\|=1$, if we denote $\eU_{n,{\cal A}} \equiv n^{-1} \sum^n_{i=1} \eeX_{i,{\cal A}}  \eeX_{i,{\cal A}}^t$, then,  $\sqrt{n} (\eu^t \eU^{-1}_{n,{\cal A}} \eu)^{-1/2} \eu^t ( \widehat{\eb}_n - \eb^0)_{\cal{A}}  \overset{\cal L} {\underset{n \rightarrow \infty}{\longrightarrow}} {\cal N}\big(0, \tau (1- \tau ) f^{-2}(0) \big)$.
\end{theorem}

We observe that, in respect to the case $p$ fixed, now we first prove the sparsity property. For showing  $\PP [   {\cal A}  \subseteq \widehat{{\cal A}}_n ] \rightarrow 1$, we prove that:
$
\lim_{n \rightarrow \infty}\PP \big[   \min_{j \in {\cal A}} \| \widehat{\eb}_{n;j} \| > 0  \big] = 1$.
For showing  $\lim_{n \rightarrow \infty}\PP [Card(  {\cal A}^c \cap \widehat{{\cal A}}_n) \geq 1]=0 $, we  use the asymptotic properties of quantile process  and imposed conditions for the tuning parameters. In order to proof the asymptotic normality of  $(\widehat{\eb}_n )_{\cal{A}}$, we mainly use  the sparsity property and we prove that for the penalized process $L_n(\eb)$, with $\eb$ in a $n^{-1/2}$-neighbourhood of $\ebo$, the penalties are much smaller than the quantile process. Finally, a CLT for the independent random variable sequences is applied.

\begin{remark}
Results of Lemma \ref{th_vconv} and of Theorem \ref{Theorem 2SPL} are new even for the particular case of quantile model without grouped variables.
\end{remark}
\begin{remark}
Algorithm and the related numerical part are a very difficult task, firstly since in the process $G_n(\eb)$ and in the two penalties of (\ref{eq3}), the variables  are grouped. On the other hand, quantile process and penalties are continuous but not differentiable in respect to parameters $\eb$. The author has not found any numerical work, even for the particular case $d=1$, of ungrouped variables, for a linear quantile model, with adaptive fused LASSO penalty. Consequently, for the  method  proposed in the present paper, another work should be conducted on numerical method, firstly for a quantile model without grouped variables and afterwards for quantile model with grouped variables. 
\end{remark}

\section{Proofs}
\label{preuves}
\hh In this section, the proofs of Lemmas and of Theorems presented is Sections \ref{section 3} and \ref{section 4} are presented.\\
In order to study the asymptotic properties of the \textit{afg\_LASSO\_Q} estimator $\widehat{\eb}_n$, we consider the following random variable
\begin{equation} 
\label{Di}
 {\cal D}_i  \equiv  (1-\tau) \e1_{\varepsilon_i <0}- \tau \e1_{\varepsilon_i \geq 0}.
  \end{equation}
 Obviously, $\eE[{\cal D}_i]=0$ and $\rho_\tau(\varepsilon_i)=- \varepsilon_i {\cal D}_i$. 
\subsection{Result proofs for   c=0 case}
\label{subsection 5.1}
\hh We start be giving the proofs of  results presented in Section \ref{section 3}.\\

\noindent {\bf Proof of Lemma \ref{Lemma 2.1}}. 
 We show that for all $\epsilon > 0$, there exists a constant  $B_\epsilon >0$ (without loss of generality, we take $B_\epsilon >0$, otherwise we take   $|B_\epsilon|$) sufficiently large such that for $ n $ large enough:
\begin{equation}
\label{eq8}
\PP \cro{ \inf_{\| \eu \| =1} L_n\pth{\eb^0+B_\epsilon n^{-1/2} \eu}>0} \geq 1-\epsilon ,
\end{equation}
with $\eu \in \R^r$, $\| \eu \| =1$.\\
Let $C_1>0$ be some constant. We will study the random process: $L_n\pth{\eb^0+C_1 n^{-1/2} \eu}=G_n(\ebo+C_1 n^{-1/2} \eu)-G_n(\ebo)+\mu_n^{(1)} \sum^p_{j=1} \|\widetilde{\eb}_{n;j} \|^{- \gamma} \big[\|\eb^0_j+n^{-1/2}C_1\eu_j \|- \|\eb^0_j\|\big]+\mu_n^{(2)} \sum^p_{j=2} \widehat{\omega}_{n;j}^{(2)} \big( \| \eb^0_j-\eb^0_{j-1} +C_1 n^{-1/2}(\eu_j-\eu_{j-1}) \|$ $ - \| \eb^0_j-\eb^0_{j-1}\|\big)$.\\
 For each observation $i$, consider   the random  variable  ${\cal R}_i\equiv \rho_\tau(\varepsilon_i -C_1 n^{-1/2} \eeX^t_i \eu) - C_1 n^{-1/2} {\cal D}_i \eeX^t_i \eu$, with ${\cal D}_i$ defined by (\ref{Di}). Consider also the following   random vector  $\textbf{W}_n \equiv C_1 n^{-1/2} \sum^n_{i=1} {\cal D}_i \eeX^t_i $. 
Then the loss term of the random process  $L_n\pth{\eb^0+C_1 n^{-1/2} \eu}$ can be written:
\begin{equation}
\label{Ge}
G_n(\ebo+C_1 n^{-1/2} \eu)-G_n(\ebo)=\eE \cro{ G_n(\ebo+C_1 n^{-1/2} \eu)-G_n(\ebo)}+\textbf{W}_n \eu+\sum^n_{i=1}({\cal R}_i- \eE[{\cal R}_i]).
\end{equation}
For the first term of the right-hand side of (\ref{Ge}) we have: 
\[
\eE \cro{ G_n\pth{\eb^0+C_1  \frac{\eu}{\sqrt{n}} } - G_n(\eb^0) }=\sum^n_{i=1} \int^{C_1 n^{-1/2} \eeX^t_i \eu }_0 [F(t)-F(0)] dt.
\]
Since $\| \eu\|=1$, by assumption (A3), we have, $n^{-1/2} \eeX^t_i \eu =o(1)$. Using assumption (A1) together with relation (\ref{rA2}), by the mean value theorem, we obtain:
\begin{equation}
\label{eq9}
 \eE \cro{ G_n\pth{\eb^0+C_1  \frac{\eu}{\sqrt{n}}} - G_n(\eb^0) } = C^2_1 \frac{f(0)}{2}   \frac{1}{n} \sum^n_{i=1} (\eeX^t_i \eu)^2 (1+o(1)).
\end{equation}
For the third term of the right-hand side of (\ref{Ge}), since the errors $\varepsilon_i$ are i.i.d., we have,
\begin{equation}
\label{ee1}
\begin{array}{c}
\displaystyle{\eE \cro{\sum^n_{i=1}({\cal R}_i- \eE[{\cal R}_i])}^2 \leq \sum^n_{i=1} \eE [{\cal R}_i^2] \leq \sum^n_{i=1}  \eE \cro{\pth{C_1 n^{-1/2} |\eeX^t_i \eu| \e1_{|\varepsilon_i| < C_1 n^{-1/2} |\eeX^t_i \eu|}}^2}} \\
 \displaystyle{\qquad \qquad  \leq C_1^2  n^{-1} \sum^n_{i=1}|\eeX^t_i \eu|^2 \eE \cro{\e1_{|\varepsilon_i| < C_1 n^{-1/2} |\eeX^t_i \eu|}}.}
 \end{array}
\end{equation}
But, using assumption (A3), 
\begin{equation}
\label{ee2}
 \eE \cro{\e1_{|\varepsilon_i| < C_1 n^{-1/2} |\eeX^t_i \eu|}} \leq C n^{-1/2} \| \eeX_i\| \leq C n^{-1/2} \max_{1 \leqslant i \leqslant n} \| \eeX_i\| =o(1).
\end{equation}
Using assumption (A2),   relations (\ref{ee1}) and (\ref{ee2})  imply: $\eE \cro{\sum^n_{i=1}({\cal R}_i- \eE[{\cal R}_i])}^2 \leq o(1)$. Then, by Bienaymé-Tchebychev inequality, we have
\begin{equation}
\label{Ro}
\sum^n_{i=1}({\cal R}_i- \eE[{\cal R}_i]) =o_{\PP}(1).
\end{equation}
 For the second term of the right-hand side of (\ref{Ge}) we have that random variable  $\textbf{W}_n \eu$ converges in distribution to a centred Gaussian law. Then, taking also into account relations (\ref{eq9}) and (\ref{Ro}), we obtain that  relation (\ref{Ge}) becomes:
\begin{equation}
\label{Gebis}
G_n(\ebo+C_1 n^{-1/2} \eu)-G_n(\ebo)=\pth{C_1^2 \frac{f(0)}{2} \frac{1}{n} \sum^n_{i=1} (\eeX^t_i \eu)^2 }(1+o_{\PP}(1)).
\end{equation}

Now we study the  penalty terms for $L_n\pth{\eb^0+C_1 n^{-1/2} \eu}$.
\begin{itemize}
\item For the penalty $\mu_n^{(1)} \sum^p_{j=1} \|\widetilde{\eb}_{n;j} \|^{- \gamma} \cro{\|\eb^0_j+n^{-1/2}C_1\eu_j \|- \|\eb^0_j\|}$,  two cases are considered for the index $j$:
\begin{itemize}
 \item if $j \in {\cal A}$, then,  since the quantile estimator $\widetilde{\eb}_{n;j}$ is consistent, we have with probability converging to 1 as  $n \rightarrow \infty$, that $\mu_n^{(1)} \left| \|\widetilde{\eb}_{n;j} \|^{- \gamma} \cro{\|\eb^0_j+n^{-1/2}C_1 \eu_j \|- \|\eb^0_j\|} \right| < C\mu_n^{(1)} n^{-1/2} \| \eu_j \|  \rightarrow 0$, by conditions of (\ref{cond_lambda}).\\
 \item if $j \in {\cal A}^c$, then, taking into account that the convergence rate of $\widetilde{\eb}_{n;j} $ to $\textbf{0}$ is $n^{-1/2}$,  this penalty is $O_{\PP}\pth{\mu_n^{(1)} \|\widetilde{\eb}_{n;j} \|^{- \gamma}  n^{-1/2}\|\eu_j \|} = O_{\PP}\pth{ \mu_n^{(1)} n^{(\gamma -1)/2}\|\eu_j \| }$,  which converges to $\infty $ when $\|\eu_j \| \neq 0$  by (\ref{cond_lambda}) and it is equal to  $0$ when $\|\eu_j \| = 0$.
 \end{itemize}
 \item We will now study the penalty $\mu_n^{(2)} \sum^p_{j=2} \widehat{\omega}_{n;j}^{(2)} \big( \| \eb^0_j-\eb^0_{j-1} +C_1 n^{-1/2}(\eu_j-\eu_{j-1}) \|$ $ - \| \eb^0_j-\eb^0_{j-1}\|\big)$. We consider the two possible cases for the index  $j$:
 \begin{itemize}
  \item if $\eb^0_j = \eb^0_{j-1}$, then we have, $\mu_n^{(2)} \widehat{\omega}_{n;j}^{(2)} \big( \| \eb^0_j-\eb^0_{j-1} +C_1n^{-1/2}(\eu_j-\eu_{j-1})\| - \| \eb^0_j-\eb^0_{j-1}\|\big)$ $=O_{\PP}\big(\mu_n^{(2)} \widehat{\omega}_{n;j}^{(2)} n^{-1/2}$ ${\| \eu_j-\eu_{j-1} \|}\big)>0$.\\
  \item if $\eb^0_j \neq \eb^0_{j-1}$, then, using conditions (\ref{cond_lambda}), we obtain, $\mu_n^{(2)} \widehat{\omega}_{n;j}^{(2)} \big( \| \eb^0_j-\eb^0_{j-1} +C_1 n^{-1/2} (\eu_j-\eu_{j-1}) \| - \| \eb^0_j-\eb^0_{j-1}\|\big)$  $=O_{\PP}\big(\mu_n^{(2)} \widehat{\omega}_{n;j}^{(2)} (\eu_j -\eu_{j-1})^t n^{-1/2}(\eb^0_j-\eb^0_{j-1} ) \| \eb^0_j-\eb^0_{j-1}\|^{-1}\big)=O_{\PP}\big(n^{-1/2}\mu_n^{(2)}( \eb^0_j-\eb^0_{j-1} )\big) =o_{\PP}(1)$.
  \end{itemize}
 \end{itemize}
 
 Then, since in the following relation  $O_{\PP}(\mu_n^{(m)} n^{(\gamma -1)/2}) >0$, for any $ m=1, 2$, and taking into account relation (\ref{Gebis}) together with   the study realised on the penalties, we have for $n$ and  $B_\epsilon $ large enough that:
 \[
 \begin{array}{c}
\displaystyle{ L_n\pth{\eb^0+B_\epsilon n^{-1/2} \eu}= B_\epsilon^2 \frac{f(0)}{2n}  \sum^n_{i=1} (\eeX^t_i \eu)^2+B_\epsilon \bigg(O_{\PP}(\mu_n^{(1)} n^{(\gamma -1)/2}) +O_{\PP}(\mu_n^{(1)} n^{-1/2}) }\\
\qquad \displaystyle{+O_{\PP}(\mu_n^{(2)} n^{-1/2})+ O_{\PP}(\mu_n^{(2)} n^{(\gamma -1)/2})\bigg) .}
 \end{array}
  \]
 Taking into account (\ref{cond_lambda}), we obtain relation (\ref{eq8}) for $n$ and $B_\epsilon$  large enough.
\hspace*{\fill}$\blacksquare$  \\

\noindent {\bf Proof of Theorem \ref{theorem 1}}. For $\eu \in \R^r$, let us consider the random process:
$
L_n(\ebo+ n^{-1/2}\eu)$, with the  process $L_n$ defined by relation (\ref{Lnb}) and $\eu \in \R^r$.\\
Let's recall that  $\widehat{\eu}_n=\sqrt{n}(\widehat{\eb}_n -\eb^0)$ is the  minimizer  in $\eu$ de $L_n(\ebo+ n^{-1/2}\eu)$. In view of  the convergence rate of the estimator $\widehat{\eb}_n$ obtained by Lemma \ref{Lemma 2.1}, we will consider $\eu\equiv (\eu_1, \cdots, \eu_p) $ bounded.  On the other hand, the process $L_n(\ebo+ n^{-1/2}\eu)$ can be written:
\begin{equation}
\label{Lnb1}
\begin{array}{c}
\displaystyle{L_n(\ebo+ n^{-1/2}\eu)=[\bfz^t_n \eu+B_n(\eu)]+ \mu_n^{(1)} \sum^p_{j=1} \widehat{\omega}_{n;j}^{(1)} \cro{\|\eb^0_j+n^{-1/2}\eu_j \|- \|\eb^0_j\|} \frac{\sqrt{n}}{\sqrt{n}} }\\
\displaystyle{+\mu_n^{(2)} \sum^p_{j=2} \widehat{\omega}_{n;j}^{(2)} \pth{ \| \eb^0_j +\frac{\eu_j}{\sqrt{n}}-\pth{\eb^0_{j-1} +\frac{\eu_{j-1}}{\sqrt{n}}}\| - \| \eb^0_j-\eb^0_{j-1}\|}\frac{\sqrt{n}}{\sqrt{n}},}
\end{array}
\end{equation}
with
\[
\bfz_n  \equiv     \frac{1}{\sqrt{n}} \sum^n_{i=1} \eeX_i  {\cal D}_i,  \qquad 
B_n(\eu)  \equiv  \sum^n_{i=1} \int^{\eeX^t_i \eu/\sqrt{n}}_0 [\e1_{\varepsilon_i < t} - \e1_{\varepsilon_i < 0} ]dt,
\]
and the random variable  ${\cal D}_i$ defined by (\ref{Di}). Since $\eE[{\cal D}_i]=0$, we have that $\eE[\bfz_n]=\textbf{0}_r$. For the loss term (the first bracket of the right-hand side) of (\ref{Lnb1}), by the CLT, using (A1), (A2) and (A3), we have
 \begin{equation}
 \label{ZB}
   \bfz_n^t \eu \overset{\cal L} {\underset{n \rightarrow \infty}{\longrightarrow}} \bfz^t \eu, \qquad 
 B_n(\eu) \overset{\PP} {\underset{n \rightarrow  \infty}{\longrightarrow}}  \frac{1}{2} f(0) \eu^t \eU \eu,
 \end{equation}
 with the random $r$-vector $\bfz \sim {\cal N}(\textbf{0}_r, \tau(1-\tau)\eU )$.\\
 We now study the two penalties of the right-hand side of  (\ref{Lnb1}).\\
 For the first penalty term of  $L_n(\ebo+ n^{-1/2}\eu)$ of (\ref{Lnb1}), we have, using the conditions of relation (\ref{cond_lambda}), that,
\begin{equation}
\label{eq1112}
\mu_n\sum^p_{j=1}  \widehat{\omega}_{n;j}  \bigg( \left\|\eb^0_j +\frac{\eu_j}{\sqrt{n}} \right\| - \| \eb^0_j \|  \bigg) \frac{\sqrt{n}}{\sqrt{n}}  \overset{\PP} {\underset{n \rightarrow \infty}{\longrightarrow}} \sum^p_{j=1}W^{(1)}(\eb^0_j,\eu),
\end{equation}
with
 \[
 W^{(1)}(\eb^0_j;\eu_j)\equiv \left\{
 \begin{array}{lll}
 0, &\textrm {if } & \eb^0_j \neq \textbf{0}_{d} \\
 0,& \textrm {if } & \eb^0_j = \textbf{0}_{d}  \textrm{ and } \eu_j =\textbf{0}_{d}\\
 \infty , & \textrm {if } & \eb^0_j = \textbf{0}_{d}  \textrm{ and } \eu_j \neq \textbf{0}_{d}.
 \end{array}
 \right.
 \]
 For the second penalty term of  $L_n(\ebo+ n^{-1/2}\eu)$ of (\ref{Lnb1}),  consider the following notations:\\
 ${\cal P}_{2,j} \equiv \sqrt{n}  \pth{ \| \eb^0_j -\eb^0_{j-1} +n^{-1/2}(\eu_j-\eu_{j-1})\| - \| \eb^0_j-\eb^0_{j-1}\|}$ and ${\cal S}_{2,j} \equiv \mu_n^{(2)} \widehat{\omega}_{n;j}^{(2)} n^{-1/2}{\cal P}_{2,j} $.\\
 For $\eb^0_j$, $\eb^0_{j-1}$ two cases are possibles.\\
\underline{\textit{If $\eb^0_j = \eb^0_{j-1}$}}, then $ {\cal P}_{2,j}=\| \eu_j-\eu_{j-1}\| \leq C  $. On the other hand, since $\|\widetilde{\eb}_{n,j} -\widetilde{\eb}_{n,j-1}  \|=O_{\PP}(n^{-1/2})$, using conditions (\ref{cond_lambda}), we have
$
\mu_n^{(2)} \widehat{\omega}_{n;j}^{(2)}n^{-1/2}=C \mu_n^{(2)}  \|\widetilde{\eb}_{n,j} -\widetilde{\eb}_{n,j-1}  \|^{ -\gamma}n^{-1/2} $ $\overset{\PP} {\underset{n \rightarrow  \infty}{\longrightarrow}}  \infty$. 
Then
\[
{\cal S}_{2,j}   \overset{\PP} {\underset{n \rightarrow  \infty}{\longrightarrow}} 
\left\{
\begin{array}{lll}
 0, &\textrm {if } & \eu_j = \eu_{j-1} \\
  \infty ,  &\textrm {if } & \eu_j \neq \eu_{j-1} .
 \end{array}
\right.
\]
\underline{\textit{If $\eb^0_j \neq \eb^0_{j-1}$}}, since
$\lim_{n \rightarrow \infty}
{\cal P}_{2,j} =   (\eu_j-\eu_{j-1})^t(\eb^0_j -\eb^0_{j-1}) \| \eb^0_j -\eb^0_{j-1}\|^{-1}$,
we have that,  ${\cal S}_{2,j} \overset{\PP} {\underset{n \rightarrow  \infty}{\longrightarrow}}  0$.\\
So, considering both cases, we can write
\begin{equation}
\label{S2j}
\sum^p_{j=1} {\cal S}_{2,j} \overset{\PP} {\underset{n \rightarrow  \infty}{\longrightarrow}}  W^{(2)}(\eb^0_j, \eb^0_{j-1};\eu_j,\eu_{j-1}),  
\end{equation}
with
\[
 W^{(2)}(\eb^0_j, \eb^0_{j-1};\eu_j,\eu_{j-1}) \equiv \left\{
 \begin{array}{lll}
 0, &\textrm {if } & \eb^0_j \neq \eb^0_{j-1} \\
 0,& \textrm {if } & \eb^0_j = \eb^0_{j-1}  \textrm{ and } \eu_j =\eu_{j-1}\\
 \infty , & \textrm {if } & \eb^0_j = \eb^0_{j-1}  \textrm{ and } \eu_j \neq \eu_{j-1}.
 \end{array}
 \right.
\]
Thus, for process (\ref{Lnb1}), taking into account of relations  (\ref{ZB}), (\ref{eq1112}) and (\ref{S2j}), we obtain,
\[
L_n(\ebo+ n^{-1/2}\eu) \overset{\cal L} {\underset{n \rightarrow \infty}{\longrightarrow}}  L(\eu), 
\]
with, the limit random variable,
\[
L(\eu) \equiv \bfz^t \eu+\frac{1}{2} f(0) \eu^t \eU \eu +\sum^p_{j=1} W^{(1)}(\eb^0_j;\eu_j)+\sum^p_{j=2}  W^{(2)}(\eb^0_j, \eb^0_{j-1};\eu_j,\eu_{j-1}).
\]
But $\widehat{\eu}_n =\argmin_{\eu} L_n(\eu)$. On the other hand, $L(\eu)$ is bounded  for any $j$ such that $\eb^0_j \neq \textbf{0}_{d}$,  when
$\eb^0_j \neq \eb^0_{j-1}$ or when $( \eb^0_j = \eb^0_{j-1}, \eu_j =\eu_{j-1})$. In these cases, the expression of $L(\eu)$ is: $L(\eu) \equiv \bfz^t \eu+ 2^{-1}f(0) \eu^t \eU \eu$.  Since the minimizer of $L(\eu)$ is the gaussian vector $f^{-1}(0) \bfz$, we deduct that  $\widehat{\eu}_n $ is  asymptotically Normal and the theorem follows.

\hspace*{\fill}$\blacksquare$ \\

\noindent {\bf Proof of Theorem \ref{th_selection}}.
By Theorem \ref{theorem 1}, for any  $  j \in {\cal A}$ we have that $\sqrt{n} (\widehat{\eb}_{n;j} - \eb^0_j) \overset{\cal L} {\underset{n \rightarrow \infty}{\longrightarrow}} {\cal N}(\textbf{0}_{d}, \tau(1- \tau) f^{-2}(0) \eU_{{\cal A}_j} ) $, with $\eU_{{\cal A}_j}$ a square matrix of size $d \times d$, the submatrix of $\eU$. Since $\eb^0_j \neq \textbf{0}_{d}$, then $j \in \widehat{{\cal A}}_n$. Thus
\begin{equation}
\label{ecA1}
\lim_{n \rightarrow \infty}\PP [{\cal A} \subseteq \widehat{{\cal A}}_n] =  1 .
\end{equation}
To finish the proof, we show that, $\lim_{n \rightarrow \infty}\PP [Card(  {\cal A}^c \cap \widehat{{\cal A}}_n) \geq 1]=0 $. We assume without loss of generality that ${\cal A}^c \cap \widehat{{\cal A}}_n=\{ j_1\}$. If this intersection contains more than one element,  the calculations are the same, except they are painful.\\
In addition to the estimator $\widehat{\eb}_n$ which   has the  $j_1$th group such that ${\cal A}^c \cap \widehat{{\cal A}}_n=\{ j_1\}$, let us consider a second estimator  $\eb^*$ for $\eb$. Taking $\eb^* \equiv (\widehat{\eb}_{n;{\cal A}}, \textbf{0}_{{\cal A}^c})$,  we will show that $Q_n(\widehat{\eb}_n) > Q_n(\eb^*)$ with a probability converging to 1.\\
For this, we will study the following difference:
\[
Q_n(\widehat{\eb}_n) - Q_n(\eb^*)=\sum^n_{i=1}  \bigg[ \rho_\tau(Y_i-  \eeX^t_i \widehat{\eb}_n) -\rho_\tau(\varepsilon_i) \bigg] -\sum^n_{i=1}  \bigg[ \rho_\tau(Y_i-  \eeX^t_i  \eb^*) -\rho_\tau(\varepsilon_i) \bigg]+\mu_n^{(1)} \widehat{\omega}_{n;j_1}^{(1)} \|\widehat{\eb}_{n;j_1} \|
\] 
\[
 +\mu_n^{(2)} \bigg[ \widehat{\omega}_{n;j_1}^{(2)} \|\widehat{\eb}_{n;j_1}- \widehat{\eb}_{n;j_1-1} \| +\widehat{\omega}_{n;j_1+1}^{(2)} \|\widehat{\eb}_{n;j_1+1}- \widehat{\eb}_{n;j_1} \|-  \widehat{\omega}_{n;j_1}^{(2)}\|\widehat{\eb}_{n;j_1-1}\| -\widehat{\omega}_{n;j_1+1}^{(2)}\|\widehat{\eb}_{n;j_1+1} \| \bigg].
\]
Similarly as in the proof of Theorem  \ref{theorem 1}, we have that, $\sum^n_{i=1}  \bigg[ \rho_\tau(Y_i-  \eeX^t_i \widehat{\eb}_n) -\rho_\tau(\varepsilon_i) \bigg]$ and $\sum^n_{i=1}  \bigg[ \rho_\tau(Y_i-  \eeX^t_i  \eb^*) -\rho_\tau(\varepsilon_i) \bigg] $ are bounded, with a probability converging to 1. Then,  with a probability converging to 1, we have that, 
\begin{equation}
\label{QQ}
\begin{array}{c}
\displaystyle{Q_n(\widehat{\eb}_n) - Q_n(\eb^*)= C +\mu_n^{(1)} \widehat{\omega}_{n;j_1}^{(1)} \|\widehat{\eb}_{n;j_1} \| +\mu_n^{(2)} \bigg[ \widehat{\omega}_{n;j_1}^{(2)}  \big(\|\widehat{\eb}_{n;j_1}- \widehat{\eb}_{n;j_1-1} \| - \|\widehat{\eb}_{n;j_1-1}\| \big)}\\
\displaystyle{+\widehat{\omega}_{n;j_1+1}^{(2)} \big(\|\widehat{\eb}_{n;j_1+1}- \widehat{\eb}_{n;j_1} \|- \|\widehat{\eb}_{n;j_1+1} \| \big) \bigg].}
\end{array}
\end{equation}
Since $j_1 \in {\cal A}^c$, we have that $\widetilde{\eb}_{n;j_1} \overset{\PP} {\underset{n \rightarrow \infty}{\longrightarrow}} 0$ and then $\widehat{\omega}_{n;j_1}^{(1)}  \overset{\PP} {\underset{n \rightarrow \infty}{\longrightarrow}} \infty$. For the first term of the penalty of  (\ref{QQ}), by conditions (\ref{cond_lambda}), we have: 
 $\mu_n^{(1)} \widehat{\omega}_{n;j_1}^{(1)} \|\widehat{\eb}_{n;j_1} \|
= \mu_n^{(1)} n^{-1/2} \|n^{1/2}\widehat{\eb}_{n;j_1} \| n^{\gamma/2}$ $ \cdot \|n^{1/2}  \widetilde{\eb}_{n;j_1}\|^{-\gamma} =O_{\PP} \pth{\mu_n^{(1)} n^{(\gamma-1)/2}} \overset{\PP} {\underset{n \rightarrow \infty}{\longrightarrow}} \infty$. \\

We will study now the two penalties, in $\mu_n^{(1)}$ and in $\mu_n^{(2)}$ of (\ref{QQ}). For $\eb^0_{j_1-1}$, $\eb^0_{j_1}$, $\eb^0_{j_1+1}$ three cases are possibles.\\
\underline{\textit{Case 1.}} $\eb^0_{j_1} \neq \eb^0_{j_1-1}$ and $\eb^0_{j_1+1} \neq \eb^0_{j_1}$.\\
In this  case,  we have that, the weights  $\widehat{\omega}_{n;j_1}^{(2)}$ and $\widehat{\omega}_{n;j_1+1}^{(2)}$ converge in probability  for $n \rightarrow\infty$ to a strictly positive bounded constant.  On the other hand, we have the following obvious inequalities, with probability 1: $\|\widehat{\eb}_{n;j_1}- \widehat{\eb}_{n;j_1-1} \|- \| \widehat{\eb}_{n;j_1-1}\| \geq - \|\widehat{\eb}_{n;j_1} \| $ and $\|\widehat{\eb}_{n;j_1+1}- \widehat{\eb}_{n;j_1} \|-\|\widehat{\eb}_{n;j_1+1}\| \geq - \|\widehat{\eb}_{n;j_1} \| $. 
Then, with probability 1, we have,
\[
\mu_n^{(1)} \widehat{\omega}_{n;j_1}^{(1)} \|\widehat{\eb}_{n;j_1} \| +\mu_n^{(2)}   \widehat{\omega}_{n;j_1}^{(2)}  \big(\|\widehat{\eb}_{n;j_1}- \widehat{\eb}_{n;j_1-1} \| - \|\widehat{\eb}_{n;j_1-1}\| \big) \geq \cro{\mu_n^{(1)} \widehat{\omega}_{n;j_1}^{(1)} -\mu_n^{(2)}  \widehat{\omega}_{n;j_1}^{(2)} }\|\widehat{\eb}_{n;j_1} \|
\]
and since $\|n^{1/2}\widehat{\eb}_{n;j_1} \|$, $\|n^{1/2} \widetilde{\eb}_{n; j_1} \|^{-\gamma}$ are bounded with probability converging to 1, since $n^{(\gamma-1)/2}\mu_n^{(1)}  $ $ \rightarrow \infty$,   we have that the above  relation is
\[
=n^{-1/2} \big[\mu_n^{(1)} \widehat{\omega}_{n;j_1}^{(1)} -\mu_n^{(2)}  \widehat{\omega}_{n;j_1}^{(2)} \big]\|n^{1/2}\widehat{\eb}_{n;j_1} \| =O_{\PP}\big( n^{-1/2}\mu_n^{(1)} \widehat{\omega}_{n;j_1}^{(1)}\big)
\]
\[
= O_{\PP}\big(n^{-1/2}\mu_n^{(1)} n^{\gamma /2} \|n^{1/2} \widetilde{\eb}_{n; j_1} \|^{-\gamma}\big) \overset{\PP} {\underset{n \rightarrow \infty}{\longrightarrow}} \infty .
\]
For the following  term of the  penalty,  we have:
\[
\begin{array}{c}
\displaystyle{ \mu_n^{(2)} \widehat{\omega}_{n;j_1+1}^{(2)} \big(\|\widehat{\eb}_{n;j_1+1}- \widehat{\eb}_{n;j_1} \|- \|\widehat{\eb}_{n;j_1+1} \| \big) \geq - \mu_n^{(2)} \widehat{\omega}_{n;j_1+1}^{(2)}\|\widehat{\eb}_{n;j_1} \| \qquad \qquad\qquad \qquad \qquad  \qquad \qquad } 
\\
 \displaystyle{\qquad \qquad \qquad \qquad \qquad \qquad \qquad  =O_{\PP}\big(-C n^{-1/2}\mu_n^{(2)}\|n^{1/2}\widehat{\eb}_{n;j_1} \| \big)  =o_{\PP}(1).}
\end{array}
\]
We have used the fact that $\widehat{\omega}_{n;j_1+1}^{(2)}$ converges in probability, for  $n \rightarrow\infty$, to a strictly positive bounded constant, $n^{-1/2}\mu_n^{(2)}\rightarrow 0$ and $\|n^{1/2}\widehat{\eb}_{n;j_1} \|$ is bounded, with probability converging to 1.\\
Thus, for relation (\ref{QQ}),  since $\mu_n^{(1)}  n^{(\gamma-1)/2}\|n^{1/2} \widetilde{\eb}_{n; j_1} \|^{-\gamma} \overset{\PP} {\underset{n \rightarrow \infty}{\longrightarrow}} \infty$, we have with a probability converging to 1, that
\begin{equation}
\label{C1}
Q_n(\widehat{\eb}_n) - Q_n(\eb^*) > C +\mu_n^{(1)}  n^{(\gamma-1)/2}\|n^{1/2} \widetilde{\eb}_{n; j_1} \|^{-\gamma}  +o_{\PP}(1) >0.
\end{equation}
\underline{\textit{Case 2.}} $\eb^0_{j_1} \neq \eb^0_{j_1-1}$ and $\eb^0_{j_1+1} = \eb^0_{j_1}$ (or vice versa, the calculations are the same). \\
Since $\eb^0_{j_1}=\textbf{0}$, then $\eb^0_{j_1+1} = \textbf{0}$ and $\eb^0_{j_1-1} \neq \textbf{0}$. \\
Since ${\cal A}^c \cap \widehat{{\cal A}_n} =\{j_1\}$, then $j_1+1 \not \in \widehat{{\cal A}_n}$ and thus $\widehat{\eb}_{n;j_1+1}=\textbf{0}$. Therefore, the weight $\widehat{\omega}_{n;j_1}^{(2)}$ converges for $n \rightarrow\infty$ to a strictly positive bounded constant and  $\widehat{\omega}_{n;j_1+1}^{(2)} \overset{\PP} {\underset{n \rightarrow \infty}{\longrightarrow}} \infty$. Then
\[
\mu_n^{(1)} \widehat{\omega}_{n;j_1}^{(1)} \|\widehat{\eb}_{n;j_1} \| +\mu_n^{(2)}  \big[ \widehat{\omega}_{n;j_1}^{(2)}  \big(\|\widehat{\eb}_{n;j_1}- \widehat{\eb}_{n;j_1-1} \| - \|\widehat{\eb}_{n;j_1-1}\| \big)+   \widehat{\omega}_{n;j_1+1}^{(2)}  \|\widehat{\eb}_{n;j_1} \|      \big]
\]
\[
\qquad \geq \mu_n^{(1)} \widehat{\omega}_{n;j_1}^{(1)} \|\widehat{\eb}_{n;j_1} \| +\mu_n^{(2)}  \bigg[- \widehat{\omega}_{n;j_1}^{(2)} \|\widehat{\eb}_{n;j_1} \|  +   \widehat{\omega}_{n;j_1+1}^{(2)}  \|\widehat{\eb}_{n;j_1} \|  \bigg] \]
\[
\qquad  \geq \mu_n^{(1)} \widehat{\omega}_{n;j_1}^{(1)} \|\widehat{\eb}_{n;j_1} \| +\mu_n^{(2)}  \widehat{\omega}_{n;j_1+1}^{(2)}  \|\widehat{\eb}_{n;j_1} \| (1+o_{\PP}(1)) \overset{\PP} {\underset{n \rightarrow \infty}{\longrightarrow}} \infty,
\]
and thus, for all $\epsilon >0$, 
\begin{equation}
\label{Cas2}
\lim_{n \rightarrow \infty} \PP\cro{Q_n(\widehat{\eb}_n) - Q_n(\eb^*) > \epsilon}=1.
\end{equation}
\underline{\textit{Case 3.}} $\eb^0_{j_1}= \eb^0_{j_1-1}=\eb^0_{j_1+1}  $.\\
In this case, we have $\widehat{\eb}_{n;j_1}  \neq \textbf{0}$, $\widehat{\eb}_{n;j_1-1} =\widehat{\eb}_{n;j_1+1}=\textbf{0}$. The penalty in this case is:
\[
\mu_n^{(1)} \widehat{\omega}_{n;j_1}^{(1)} \|\widehat{\eb}_{n;j_1} \| +\mu_n^{(2)} \bigg[\widehat{\omega}_{n;j_1}^{(2)} +\widehat{\omega}_{n;j_1+1}^{(2)} \bigg] \|\widehat{\eb}_{n;j_1}\| \overset{\PP} {\underset{n \rightarrow \infty}{\longrightarrow}}  \infty,
\]
which means for relation (\ref{QQ}), that for all $\epsilon >0$, 
\begin{equation}
\label{Cas3}
\lim_{n \rightarrow \infty} \PP\cro{Q_n(\widehat{\eb}_n) - Q_n(\eb^*) > \epsilon}=1.
\end{equation}
\\
Thus, in all three cases, taking into account (\ref{QQ}), (\ref{C1}), (\ref{Cas2}), (\ref{Cas3}), we get:
\[
\lim_{n \rightarrow \infty} \PP\cro{ Q_n(\widehat{\eb}_n) > Q_n(\eb^*)}=1,
\]
which implies that,  $\widehat{\eb}_n$ is not the minimizer of $Q_n$. Therefore, $\lim_{n \rightarrow \infty}\PP [j_1 \in {\cal A}^c \cap \widehat{{\cal A}}_n]=0 $ and the proof is finished.
\hspace*{\fill}$\blacksquare$ \\

\subsection{Result proofs for  $c>0$ case}
\label{subsection 5.2}
\hh We present now the proofs of the results stated in Section \ref{section 4}. For some results, we will consider obtained results in \cite{Ciuperca-16} for asymptotic behaviour of the non fused adaptive group LASSO quantile estimator. \\

\noindent {\bf Proof of Lemma \ref{th_vconv}}. Let be a positive constant $C_1$and a $r$-vector $\eu$ such that $\|\eu\|=1$. For the process  $L_n$ defined by  (\ref{Lnb}), let us consider
\[
L_n \pth{ \ebo+ C_1 \sqrt{\frac{p}{n}} \eu}  =
   G_n\pth{ \ebo+ C_1 \sqrt{\frac{p}{n}} \eu} -   G_n(\ebo) +\sum^p_{j=1} \mu^{(1)}_n \widehat{\omega}^{(1)}_{n;j} \cro{\|\eb^0_j+\sqrt{\frac{p}{n}} C_1 \eu_j \|- \|\eb^0_j\| }
\]
\begin{equation}
\label{T}
+\mu^{(2)}_n\sum^p_{j=2}\widehat{\omega}^{(2)}_{n;j} \cro{\|\eb^0_j+\sqrt{\frac{p}{n}} C_1 \eu_j -\pth{\eb^0_{j-1}+\sqrt{\frac{p}{n}} C_1 \eu_{j-1}} \| -\|\eb^0_j - \eb^0_{j-1}\|} \equiv T_1+T_2+T_3.
\end{equation}
For $T_1, T_2$,  using  $\mu^{(1)}_n n^{(c-1)/{2} - \alpha \gamma} \rightarrow 0$, proceeding as in the proof of   Theorem 3 of \cite{Ciuperca-16}, we get:
\begin{equation}
\label{T12}
T_1+T_2 > C_1^2 f(0)  p  \pth{\frac{1}{n}  \sum^n_{i=1} \eu^t \eeX_i \eeX_i^t \eu } (1+o_{\PP}(1)) - C_1O_{\PP}\pth{  {p} }.
\end{equation}
For $T_3$, we have:
\[
T_3  \geq \mu^{(2)}_n\sum_{\{j; \eb^0_j \neq \eb^0_{j-1}\}}\widehat{\omega}^{(2)}_{n;j}\cro{\|\eb^0_j+\sqrt{\frac{p}{n}} C_1 \eu_j -\pth{\eb^0_{j-1}+\sqrt{\frac{p}{n}} C_1 \eu_{j-1}} \| -\|\eb^0_j - \eb^0_{j-1}\|}
\]
\begin{equation}
\label{T3}
\quad \geq - C_1 \sqrt{\frac{p}{n}}\mu^{(2)}_n\sum_{\{j; \eb^0_j \neq \eb^0_{j-1}\}}\widehat{\omega}^{(2)}_{n;j} \| \eu_j - \eu_{j-1}\| \geq - C_1 O_{\PP}\pth{  p}
\end{equation}
by similar reasoning to that $T_1$, using also condition $\mu^{(2)}_n n^{(c-1)/{2} - \alpha \gamma} \rightarrow 0$ and assumptions (A4), (A5).\\
Taking into account (\ref{T}), (\ref{T12}) and (\ref{T3}), we have for $n$ and $C_1$ large enough, that for all $\epsilon >0$:
\[
\PP \cro{ \inf_{\| \eu \| =1} L_n \pth{ \ebo+ C_1 \sqrt{\frac{p}{n}} \eu}> 0} \geq 1-\epsilon ,
\]
and the theorem follows.
\hspace*{\fill}$\blacksquare$ \\

\noindent {\bf Proof of Theorem \ref{Theorem 2SPL}}.
\textit{(i)} We show, as for the case $p$ fixed, that $\lim_{n \rightarrow \infty}\PP [Card(  {\cal A}^c \cap \widehat{{\cal A}}_n) \geq 1]=0 $. Therefore we not give some calculation details.  We assume without loss of generality that ${\cal A}^c \cap \widehat{{\cal A}}_n=\{ j_1\}$. We consider the second estimator $\eb^* \equiv (\widehat{\eb}_{n;{\cal A}}, \textbf{0}_{{\cal A}^c})$ of $\eb$ and we will show that $Q_n(\widehat{\eb}_n) > Q_n(\eb^*)$ with a probability converging to 1. We then study the following difference:
\begin{equation}
\label{EE}
\begin{array}{c}
\displaystyle{ Q_n(\widehat{\eb}_n) - Q_n(\eb^*)=\sum^n_{i=1}  \big[ \rho_\tau(Y_i-  \eeX^t_i \widehat{\eb}_n) -\rho_\tau(\varepsilon_i) \big] -\sum^n_{i=1}  \big[ \rho_\tau(Y_i-  \eeX^t_i  \eb^*) -\rho_\tau(\varepsilon_i) \big] } \\
\displaystyle{+\mu_n^{(1)} \widehat{\omega}_{n;j_1}^{(1)} \|\widehat{\eb}_{n;j_1} \| +\mu_n^{(2)} \big[ \widehat{\omega}_{n;j_1}^{(2)} \|\widehat{\eb}_{n;j_1}- \widehat{\eb}_{n;j_1-1} \| +\widehat{\omega}_{n;j_1+1}^{(2)} \|\widehat{\eb}_{n;j_1+1}- \widehat{\eb}_{n;j_1} \| }\\
\displaystyle{-  \widehat{\omega}_{n;j_1}^{(2)}\|\widehat{\eb}_{n;j_1-1}\|-\widehat{\omega}_{n;j_1+1}^{(2)}\|\widehat{\eb}_{n;j_1+1} \| \big].}
\end{array}
\end{equation}
Under assumptions (A1)-(A4), we have, as in the proof of Theorem 4 of \cite{Ciuperca-16}, that: 
\begin{equation}
\label{epr}
\sum^n_{i=1}  \bigg[ \rho_\tau(Y_i-  \eeX^t_i \widehat{\eb}_n) -\rho_\tau(\varepsilon_i) \bigg] -\sum^n_{i=1}  \bigg[ \rho_\tau(Y_i-  \eeX^t_i  \eb^*) -\rho_\tau(\varepsilon_i) \bigg]=O_{\PP} \big( n\|\widehat{\eb}_n -\eb^* \|^2\big)=O_{\PP} \pth{p} =O_{\PP}(n^c).
\end{equation}
For the first term of the penalty of (\ref{EE}), we have:
\begin{equation}
\label{ep1}
0< \mu_n^{(1)} \widehat{\omega}_{n;j_1}^{(1)} \|\widehat{\eb}_{n;j_1} \|=\mu_n^{(1)} O_{\PP} \pth{  \frac{p}{n}}^{(1-\gamma)/2}= \mu_n^{(1)}O_{\PP}(n^{(c-1)(1-\gamma)/2}).
\end{equation}
Since  $j_1 \in {\cal A}^c$  and also $\mu^{(1)}_n n^{\big(-c(1+\gamma)+\gamma-1 \big)/2} \rightarrow \infty$, as $n \rightarrow \infty$,we have that  term  (\ref{ep1})   dominates (\ref{epr}).\\

Now we study the penalty corresponding to $\mu^{(2)}_n$, by considering the three possible cases for $\eb^0_{j_1-1}$, $\eb^0_{j_1}$, $\eb^0_{j_1+1}$.\\
\underline{\textit{Case 1.}} $\eb^0_{j_1} \neq \eb^0_{j_1-1}$ and $\eb^0_{j_1+1} \neq \eb^0_{j_1}$. \\
Then $\widehat{\omega}_{n;j_1}^{(2)} =\|\widetilde{\eb}_{j_1} - \widetilde{\eb}_{j_1-1} \|^{-\gamma}= \|\widetilde{\eb}_{j_1}- \eb^0_{j_1} - \widetilde{\eb}_{j_1-1}+\eb^0_{j_1-1} - \eb^0_{j_1-1}\|^{-\gamma}=O_{\PP}\big( \left\| ( {p}n^{-1})^{1/2}- \eb^0_{j_1-1}\right\|^{-\gamma}\big)$ $=O_{\PP}\big((n^{(c-1)/2}- n^\alpha)^{-\gamma}\big)= O_{\PP}\big(n^{- \alpha \gamma}\big)$. On the other hand, $\widehat{\omega}_{n;j_1}^{(1)}=O_{\PP}\big( (p n^{-1})^{-\gamma /2}\big)=O_{\PP}\big(n^{(1-c)\gamma /2}\big)$.  Since $\alpha > (c-1)/2$, then, $\widehat{\omega}_{n;j_1}^{(1)} \gg \widehat{\omega}_{n;j_1}^{(2)}$, with a probability converging to  1. Therefore, $\mu_n^{(1)} \widehat{\omega}_{n;j_1}^{(1)} \|\widehat{\eb}_{n;j_1} \| +\mu_n^{(2)}   \widehat{\omega}_{n;j_1}^{(2)}  \big(\|\widehat{\eb}_{n;j_1}- \widehat{\eb}_{n;j_1-1} \| - \|\widehat{\eb}_{n;j_1-1}\| \big) \geq \mu_n^{(1)} \widehat{\omega}_{n;j_1}^{(1)} \|\widehat{\eb}_{n;j_1} \| (1+o_{\PP}(1)) \overset{\PP} {\underset{n \rightarrow \infty}{\longrightarrow}} \infty $.\\
On the other hand, $\mu_n^{(2)} \widehat{\omega}_{n;j_1+1}^{(2)} \|\widehat{\eb}_{n;j_1} \|= \mu_n^{(2)} \|\widehat{\eb}_{n;j_1} \| \cdot \|\widetilde{\eb}_{j_1}- \widetilde{\eb}_{j_1+1} \|^{-\gamma} $, with $\|\widetilde{\eb}_{j_1}- \widetilde{\eb}_{j_1+1} \|  =O_{\PP}\bigg( ( {p}n^{-1})^{1/2}+\|-  \widetilde{\eb}_{j_1+1} + \eb^0_{j_1+1} -\eb^0_{j_1+1}\| \bigg) =O_{\PP}\pth{ n^{(c-1)/2}+ n^\alpha } =O_{\PP}\pth{ n^\alpha}$. Thus, $\widehat{\omega}_{n;j_1}^{(1)} \gg \widehat{\omega}_{n;j_1+1}^{(2)}$,  with a probability converging to  1. \\
In conclusion, the right-hand side of  $Q_n(\widehat{\eb}_n) - Q_n(\eb^*)$ of relation (\ref{EE}) is dominated by $\mu_n^{(1)} \widehat{\omega}_{n;j_1}^{(1)} \|\widehat{\eb}_{n;j_1} \|$, which converges in probability to  $\infty$.\\
\underline{\textit{Case 2.}} $\eb^0_{j_1} \neq \eb^0_{j_1-1}$ and $\eb^0_{j_1+1} = \eb^0_{j_1}=\textbf{0}$ (or vice versa). \\
In this case we have, $\widehat{\omega}_{n;j_1}^{(2)}=\| \widetilde{\eb}_{j_1}- \widetilde{\eb}_{j_1-1}\|^{-\gamma}=O_{\PP}\pth{ \|( {p}n^{-1})^{1/2}- \eb^0_{j_1-1} \|^{-  \gamma}}=O_{\PP}\pth{ n^{- \alpha \gamma} }$. Then, $ \widehat{\omega}_{n;j_1}^{(1)} \gg\widehat{\omega}_{n;j_1}^{(2)} $, with a  probability converging to  1.\\
We have also that $\widehat{\omega}_{n;j_1+1}^{(2)}=\| \widetilde{\eb}_{j_1+1}- \widetilde{\eb}_{j_1}\|^{-\gamma}=O_{\PP}\pth{ ( p n^{-1})^{- \gamma /2}}=O_{\PP}\pth{   n^{(1-c)\gamma /2}} \gg \widehat{\omega}_{n;j_1}^{(2)}$, with probability converging to  1.\\
Afterwards, relations are similar to those  of the fixed  $p$ case and we deduce (\ref{Cas2}).\\
\underline{\textit{Case 3.}}  $\eb^0_{j_1}= \eb^0_{j_1-1}=\eb^0_{j_1+1}  $, is  similar to the fixed $p$ case to derive relation (\ref{Cas3}).\\

Thus, in all three cases, we have:
$
Q_n(\widehat{\eb}_n) > Q_n(\eb^*)$  with a probability converging to one, 
which implies $\widehat{\eb}_n$ is not the minimizer of $Q_n$. Thus, $\lim_{n \rightarrow \infty}\PP [j_1 \in {\cal A}^c \cap \widehat{{\cal A}}_n]=0 $.\\

To complete the demonstration of claim \textit{(i)} we need to show that  $\PP [   {\cal A}  \subseteq \widehat{{\cal A}}_n ] \rightarrow 1$. For this, we prove that:
\[
\lim_{n \rightarrow \infty}\PP \big[   \min_{j \in {\cal A}} \| \widehat{\eb}_{n;j} \| > 0  \big] = 1.
\]
The proof is similar to that of Theorem 4 of  \cite{Ciuperca-16}.\\

\textit{(ii)} Taking into account claim \textit{(i)} and assumption (A6), the estimator $\widehat{\eb}_n$ can be written, with a probability converging to 1, as, $\widehat{\eb}_n = \ebo+ ({p}n^{-1})^{1/2} \ed$, with  $\ed=(\ed_{\cal A}, \ed_{{\cal A}^c})$, $\ed_{{\cal A}^c}=\textbf{0}_{r-r^0}$, $\| \ed_{\cal A}\| \leq C $. Then, we will consider the parameters of the  form $\eb=\ebo+({p}n^{-1})^{1/2} \ed$. \\
In order to prove claim (ii), let us consider the following penalized random process:
\begin{equation}
\label{Li}
L_n\pth{\ebo+ \sqrt{\frac{p}{n}} \ed}   =  \sum^n_{i=1} \bigg[\rho_\tau\bigg(Y_i -\eeX_i^t( \eb^0 + \sqrt{\frac{p}{n}} \ed)\bigg) - \rho_\tau(\varepsilon_i)\bigg]+{\cal P}_1 +{\cal P}_2, 
\end{equation}
with the penalties, ${\cal P}_1 \equiv  \mu_n^{(1)} \sum^p_{j=1}\widehat{\omega}_{n;j}^{(1)} \cro{\| \eb_j\| - \|\eb_{j}^0 \|}$ and  ${\cal P}_2 \equiv \mu_n^{(2)} \sum^p_{j=2} \widehat{\omega}_{n;j}^{(2)} \cro{ \| \eb_j - \eb_{j-1} \| - \| \eb_j^0 - \eb^0_{j-1} \| }$.\\
For the first term of the right-hand side of (\ref{Li}), we have
\[
 \sum^n_{i=1} \bigg[\rho_\tau\bigg(Y_i -\eeX_i^t( \eb^0 + \sqrt{\frac{p}{n}} \ed)\bigg) - \rho_\tau(\varepsilon_i)\bigg]  \qquad \qquad  \qquad \qquad \qquad \qquad
\]
\[
=  n \pth{\frac{1}{n} \sqrt{\frac{p}{n}} \sum^n_{i=1} \eeX^t_i \ed [\e1_{\varepsilon \leq 0} - \tau] + \frac{1}{n} \sum^n_{i=1} \int^{\sqrt{\frac{p}{n}} \eeX^t_i \ed}_ 0 [ \e1_{\varepsilon_i \leq t} - \e1_{\varepsilon_i \leq 0} ] dt} \equiv  n(J_1+J_2).
\]
As in the proof of Theorem 4 of \cite{Ciuperca-16}, we obtain that $\eE[J_2]=Cf(0)pn^{-1}$, also that  the difference of the  quantile processes is:
\begin{equation}
\label{Li1}
  \sum^n_{i=1} \bigg[\rho_\tau \bigg(Y_i -\eeX_i^t( \eb^0 + \sqrt{\frac{p}{n}} \ed)\bigg) - \rho_\tau(\varepsilon_i)\bigg] =   {\sqrt{\frac{p}{n}}}   \sum^n_{i=1}\eeX^t_{i,{{\cal A}}} \ed_{{\cal A}} [\e1_{\varepsilon_i <0} - \tau]+  \frac{f(0)}{2}  {p}  \ed^t_{{\cal A}} \eU_{n,{\cal A}} \ed_{{\cal A}}(1+o_{\PP}(1))
\end{equation}
and for the first penalty of  (\ref{Li}), that,
 \begin{equation}
\label{Li2}
 |{\cal P}_1| \leq O_{\PP}\pth{ \mu_n^{(1)} p n ^{(c-1)/2- \alpha \gamma} }.
 \end{equation}
 This last relation implies that  $ |{\cal P}_1| \big(n \eE[J_2]\big)^{-1}=O_{\PP}\pth{\mu_n^{(1)} n ^{(c-1)/2- \alpha \gamma} } =o_{\PP}(1)$, by condition imposed on $\mu_n^{(1)}$. \\
  Let us now consider  penalty ${\cal P}_2 $ of (\ref{Li}). Using assumptions (A4), (A5), we have, 
\[
|{\cal P}_2 | \leq C p \mu_n^{(2)} \max_{j \in {\cal A} \setminus \{ 2\}} \big(\widehat{\omega}_{n;j}^{(2)}\big) \pth{\frac{p}{n}}^{1/2} = O_{\PP}\pth{n^c n^{(c-1)/2}\mu_n^{(2)}  \max_{j \in {\cal A} \setminus \{ 2\}} \big(\widehat{\omega}_{n;j}^{(2)}\big)} = O_{\PP}\pth{n^{(3c-1)/2}\mu_n^{(2)} n^{- \alpha \gamma}}.
\]
 Thus, by condition imposed to $\mu_n^{(2)}$, we have,
\begin{equation}
\label{Li3}
\frac{|{\cal P}_2|}{n \eE[J_2]}=O_{\PP}\pth{\mu_n^{(2)} n ^{(c-1)/2- \alpha \gamma}} =o_{\PP}(1).
\end{equation}
In conclusion, taking into account relations (\ref{Li})-(\ref{Li3}), we have that:
\[
L_n\pth{\ebo+ \sqrt{\frac{p}{n}} \ed}   ={\sqrt{\frac{p}{n}}}   \sum^n_{i=1}\eeX^t_{i,{{\cal A}}} \ed_{{\cal A}} [\e1_{\varepsilon_i <0} - \tau]+  \frac{f(0)}{2}  {p}  \ed^t_{{\cal A}} \eU_{n,{\cal A}} \ed_{{\cal A}}(1+o_{\PP}(1)).
\]
The minimizer of the right-hand side of the last equation is:
\[
 \sqrt{\frac{p}{n}} \ed_{{\cal A}}=-\frac{1}{n} \frac{1}{f(0)} \eU^{-1}_{n,{\cal A}} \big( \sum^n_{i=1} \eeX_{i,{ {\cal A}}} (\e1_{\varepsilon_i \leq 0}-\tau)\big).
\]
Claim \textit{(ii)} follows by taking into account of the fact that  $\widehat \eb_{{\cal A}}- \eb^0_{{\cal A}}= ( {p}n^{-1})^{1/2} \ed_{{\cal A}}$ and by applying the CLT for the following independent random variable sequence:  
 $\pth{(f(0))^{-1} \eu^t \eU^{-1}_{n,{\cal A}}\eeX_{i {\cal A}} (\e1_{\varepsilon_i \leq 0}-\tau)}_{1 \leq i \leq n}$, with $\eu$ a $r^0$-vector such that $\|\eu\|=1$.
\hspace*{\fill}$\blacksquare$ \\


\end{document}